\documentclass[11pt]{amsart}
\usepackage{graphicx}
\usepackage{amssymb}
\usepackage{amscd}
\usepackage{cite}
\usepackage[all]{xy}
\usepackage{url}
\newtheorem*{Mtw}{Main Theorem}

\theoremstyle{remark}
\newtheorem*{uw}{Remark}

\theoremstyle{definition}

\newcommand{\Cal}[1]{\mathcal{#1}}

\newcommand{\ro}{\varrho}

\newcommand{\gen}[1]{\langle #1 \rangle}
\newcommand{\map}[3]{#1\colon #2\to #3}
\newcommand{\field}[1]{\mathbb{#1}}
\newcommand{\zz}{\field{Z}}

\newcommand{\rr}{\field{R}}
\newcommand{\cc}{\field{C}}

\newcommand{\cM}{{\Cal M}}

\begin{document}

\numberwithin{equation}{section}
\title[The hyperelliptic mapping class group of a nonorientable\ldots]
{The hyperelliptic mapping class group of a nonorientable surface of genus $g\geq 4$ has a faithful representation into  $\textrm{GL}(g^2-1,\rr)$.}

\author{Micha\l\ Stukow}

\thanks{Supported by grant 2015/17/B/ST1/03235 of National Science Centre, Poland.}
\address[]{
Institute of Mathematics, Faculty of Mathematics, Physics and Informatics, University of Gda\'nsk, 80-308 Gda\'nsk, Poland }
\email{trojkat@mat.ug.edu.pl}




\begin{abstract}
We prove that the hyperelliptic mapping class group of a nonorientable surface of genus $g\geq 4$ has a faithful linear representation of dimension $g^2-1$ over $\rr$.
\end{abstract}

\maketitle%
 \section{Introduction}%
Let $N_{g,n}$ be a smooth, nonorientable, compact surface of genus $g$ with $n$ punctures. If $n$ is zero then we omit it from the notation. Recall that $N_{g}$ is a connected sum of $g$ projective planes and $N_{g,n}$ is
obtained from $N_g$ by specifying the set $\Sigma$
of $n$ distinguished points in the interior of $N_g$.

Let ${\textrm{Diff}}(N_{g,n})$ be the group of all diffeomorphisms $\map{h}{N_{g,n}}{N_{g,n}}$ such that $h(\Sigma)=\Sigma$. By ${\Cal{M}}(N_{g,n})$ we denote the quotient group of ${\textrm{Diff}}(N_{g,n})$ by
the subgroup consisting of maps isotopic to the identity, where we assume that maps and isotopies fix the set $\Sigma$. ${\Cal{M}}(N_{g,n})$ is called the \emph{mapping class group} of $N_{g,n}$. 

The mapping class group ${\Cal{M}}(S_{g,n})$ of an orientable surface $S_{g,n}$ of genus $g$ with $n$ punctures is defined analogously, but we consider only orientation preserving maps. If we include orientation reversing maps, we obtain the so-called \emph{extended mapping class group} ${\Cal{M}}^{\pm}(S_{g,n})$.
\begin{figure}[h]
\begin{center}
\includegraphics[width=0.9\textwidth]{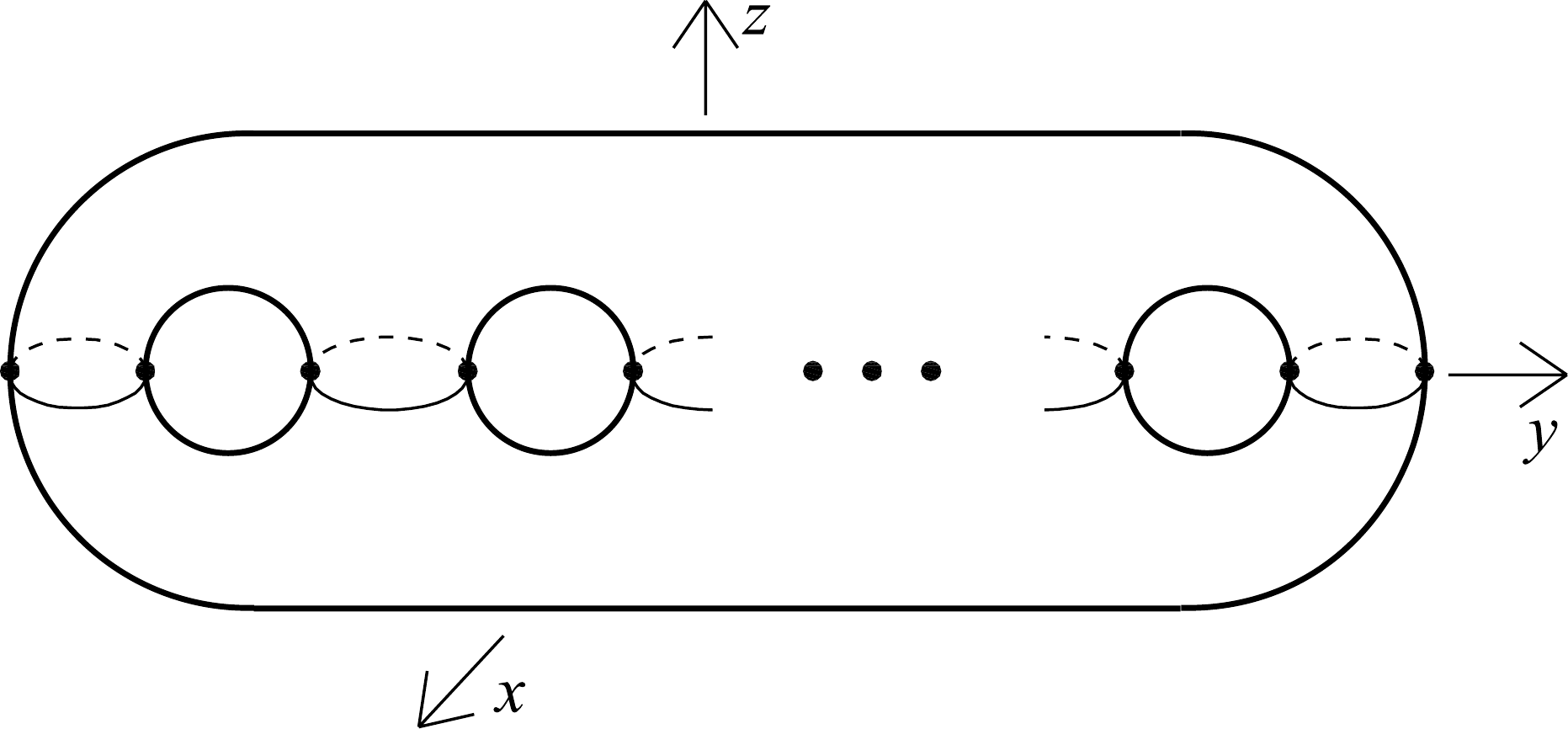}
\caption{Surface $S_g$ embedded in $\rr^3$.}\label{r01} %
\end{center}
\end{figure}
Suppose that the closed orientable surface $S_{g-1}$, where $g-1\geq 2$, is embedded in $\rr^3$ as shown in Figure \ref{r01}, in such a way that it is invariant under reflections across $xy,yz,xz$ planes. Let $\map{j}{S_{g-1}}{S_{g-1}}$ be the symmetry defined by $j(x,y,z)=(-x,-y,-z)$. Denote by $C_{\cM^{\pm}(S_{g-1})}(j)$
the centraliser of $j$ in $\cM^{\pm}(S_{g-1})$. The orbit space $S_{g-1}/\gen{j}$ is a nonorientable
surface $N_{g}$ of genus $g$ and it is known (Theorem 1 of \cite{BirChil1}) that the orbit space projection induces an epimorphism 
\[\map{\pi_j}{C_{\cM^{\pm}(S_{g-1})}(j)}{\cM(N_{g})}\]
with kernel $\ker \pi_j=\gen{j}$. In particular
\[\cM(N_{g})\cong C_{\cM^{\pm}(S_{g-1})}(j)/\gen{j}.\]
As was observed in the proof of Theorem 2.1 of \cite{Stukow_HiperOsaka}, projection $\pi_j$ has a section
 \[\map{i_j}{\cM(N_{g})}{C_{\cM(S_{g-1})}(j)\subset \cM(S_{g-1})}.\]
 In fact, for any $h\in{\cM(N_{g})}$ we can define $i_j(h)$ to be an orientation preserving lift of $h$.

Let $\ro\in C_{\cM^{\pm}(S_{g-1})}(j)$
be the \emph{hyperelliptic involution}, i.e. the half turn about the $y$-axis. The \emph{hyperelliptic mapping class group} $\cM^h(S_{g-1})$ is defined to be the centraliser of $\ro$ in $\cM(S_{g-1})$. The hyperelliptic mapping class group turns out to be a very interesting and important subgroup, in particular its finite subgroups correspond to automorphism groups of hyperelliptic Riemann surfaces -- see for example \cite{MaxHyp} and references there.

Recently we extended the notion of the hyperelliptic mapping class group to nonorientable surfaces \cite{Stukow_HiperOsaka}, by defining $\cM^h(N_g)$ to be the centraliser of $\pi_j(\ro)$ in the mapping class group $\cM(N_g)$. This definition is motivated by the notion of hyperelliptic Klein surfaces -- see for example 
\cite{Costa-hyper, Bujalance-hyper-Klein}. We say that $\pi_j(\ro)$ is the \emph{hyperelliptic
involution} of $N_g$ and by abuse of notation we write $\ro$ for $\pi_j(\ro)$.

Since $\ro\in C_{\cM^{\pm}(S_{g-1})}(j)$, we have restrictions of $\pi_j$ and $i_j$ to the maps
\[\begin{aligned}
   &\map{\pi_j}{C_{\cM^{\pm}(S_{g-1})}(\gen{j,\ro})}{\cM^h(N_{g})}\\
   &\map{i_j}{\cM^h(N_{g})}{C_{\cM(S_{g-1})}(\gen{j,\ro})\subset \cM^h(S_{g-1})}.
  \end{aligned}
\]
\section{Linear representations of the hyperelliptic mapping class group.}
Mapping class groups of projective plane $N_1$ and of Klein bottle $N_2$ are finite, hence the first nontrivial case is the group $\cM(N_3)$. This is an interesting case, because it is well known \cite{BirChil1,Scharlemann} that 
\[\cM^h(N_3)=\cM(N_3)\cong \textrm{GL}(2,\zz).\]
In particular, $\cM^h(N_3)$ has a faithful linear representation of real dimension 2. 

For $g\geq 4$, we can produce a faithful linear representation of the hyperelliptic mapping class group $\cM^h(N_g)$ as a composition of the section 
\[\map{i_j}{\cM^h(N_{g})}{C_{\cM(S_{g-1})}(\gen{j,\ro})\subset \cM^h(S_{g-1})}\]
and a faithful linear representation of $\cM^h(S_{g-1})$ obtained by Korkmaz \cite{Kork-lin} or by Bigelow and Budney \cite{Big-lin}. Recall that both of these representations of $\cM^h(S_{g-1})$ are obtained form the Lawrence--Krammer representation of the braid group \cite{BigelowLin,Krammer2}.

The above argument is immediate, but the resulting representation of $\cM^h(N_{g})$ is far from being optimal. In fact, if we use Bigelow--Budney representation of $\cM^h(S_{g-1})$ (which has much smaller dimension than the one obtained by Korkmaz) the dimension of the obtained representation of $\cM^h(N_g)$ is equal to
\[2g\cdot {2g-1\choose 2}+2(g-1)=2(g-1)(2g^2-g+1).\]
%
%
\begin{Mtw}
 If $g\geq 4$, then the hyperelliptic mapping class group $\cM^h(N_g)$ has a faithful linear representation of real dimension $g^2-1$.
\end{Mtw}
\begin{proof}
 Let $\cM^{\pm}(S_{0,g+1})$ be the extended mapping class group of a sphere with $g+1$ punctures $\{p_1,\ldots,p_{g+1}\}$, and let $\cM^{\pm}(S_{0,g,1})$ be the stabiliser of $p_{g+1}$ with respect to the action of $\cM^{\pm}(S_{0,g+1})$ on the set of punctures.
  By Theorem 2.1 of \cite{Stukow_HiperOsaka}, the orbit space projection ${\cM^h(N_g)}\to {\cM^h(N_g)}/\gen{\ro}$ induces an epimorphism 
  \[\map{\pi_\ro}{\cM^h(N_g)}{\cM^{\pm}(S_{0,g,1})}\]
with $\ker \pi_\ro=\gen{\ro}$. Moreover, by rescaling the Lawrence-Krammer representation of the braid group \cite{BigelowLin},  Bigelow and Budney constructed in the proof of Theorem 2.1 of \cite{Big-lin} a faithful linear representation 
\[\map{{\Cal{L}}'}{\cM(S_{0,g,1})}{\textrm{GL}\left({g\choose 2},\rr\right)}.\]
To be more precise, they obtained a representation over $\cc$, however their argument works without any changes over $\rr$.

Since $\cM(S_{0,g,1})$ is a subgroup of index 2 in $\cM^{\pm}(S_{0,g,1})$, the later group has an induced  faithful linear representation of dimension $2\cdot {g\choose 2}=g^2-g$.
This gives us a linear representation 
\[\map{{\Cal{L}_1}}{\cM^h(N_g)}{\textrm{GL}\left(g^2-g,\rr\right)}\]
with kernel $\ker {\Cal{L}_1}=\gen{\ro}$. It is straightforward to check that if 
\[\map{{\Cal{L}_2}}{\cM^h(N_g)}{\textrm{H}_1(N_g;\rr)}\subset {\textrm{GL}\left(g-1,\rr\right)}\]
is a standard homology representation then ${\Cal{L}_1}\oplus {\Cal{L}_2}$ is a required faithful linear representation of $\cM^h(N_g)$ of dimension $g^2-g+g-1=g^2-1$.
%
\end{proof}
\begin{uw}
The Main Theorem gives an upper bound $g^2-1$ on the minimal dimension of a faithful linear representation of the hyperelliptic mapping class group $\cM^h(N_g)$. As we mentioned in the introduction, the hyperelliptic mapping class group $\cM^h(N_3)$ has a faithful linear representation of real dimension 2, hence it seems very unlikely that the obtained bound is sharp. 
\end{uw}
\section*{Acknowledgements}
The author wishes to thank the referee for his/her helpful suggestions.

\end{document}